\documentclass[12pt,twoside,doublespace]{article}

\usepackage{latexsym,amssymb,amsfonts}
\usepackage{graphicx}

\def\smskip{\par\vskip 5 pt}
\def\QED{\hfill $\Box$\smskip}
\newtheorem{theorem}{Theorem}

\newtheorem{proposition}{Proposition}

\begin{document}

\begin{center}

\vspace{35pt}

{\Large \bf Exact Penalties for Decomposable }

\vspace{5pt}

{\Large \bf  Optimization Problems

}

\vspace{5pt}

\vspace{35pt}

{\sc Igor V.~Konnov\footnote{\normalsize E-mail: konn-igor@ya.ru}}

\vspace{35pt}

{\em  Department of System Analysis
and Information Technologies, \\ Kazan Federal University, ul.
Kremlevskaya, 18, Kazan 420008, Russia.}

\end{center}

\vspace{25pt}

\hrule

\vspace{10pt}

{\bf Abstract:} We consider a general decomposable convex optimization problem.
By using right-hand side allocation
technique, it can be transformed into a collection of small
dimensional optimization problems. The master problem is a
convex non-smooth optimization problem. We propose to apply
the exact non-smooth penalty method, which gives a solution of the initial problem
under some fixed penalty parameter
and provides the consistency of lower level problems.
The master problem is suggested to be solved by
a two-speed subgradient projection method, which
enhances the step-size selection. Preliminary results of computational experiments
confirm its efficiency.


\vspace{10pt}

\hrule

\vspace{10pt}

{\em MS Classification:} {90C25, 90C06, 90C30}


\section{Introduction} \label{s1}

The general optimization problem consists in finding the minimal
value of some goal function $f : \mathbb{R}^{n} \to \mathbb{R}$ on
a feasible set $D \subseteq \mathbb{R}^{n}$. For brevity, we write
this problem as
\begin{equation} \label{eq:1.1}
 \min \limits _{x \in D} \to f(x),
\end{equation}
its solution set is denoted by $D^{*}$ and the optimal value of the
function by $f^{*}$, i.e.
$$
f^{*} = \inf \limits _{ x \in D} f(x).
$$
It is well known that solution methods that take into account
peculiarities of particular problems show better computational
results than general purpose oriented ones. In particular, most
large dimensional optimization problems have a special structure
which admits decomposition. There exist various decomposition
approaches  allowing one to replace the initial large scale problem with a
sequence of relatively small and simple ones. Besides, they justify certain
decentralized control principles for large organization and
industrial systems; see e.g. \cite{Las70,Tsu81,GT89,CCM06} and the
references therein. Together with the known
price (Dantzig-Wolfe) and variable (Benders) decomposition
the so-called right-hand side (Kornai-Liptak)
decomposition approach was proved to be rather efficient in Optimization; see
\cite{KL65,Las70,Tsu81}. Its main idea is that
one can attain the global optimal value in a complex system by some
upper level allocation procedure of joint resource shares, whereas subsystems
(elements) are free in choosing their best actions within these shares. In particular, this
method enables us to reduce a large scale smooth
convex optimization problem to a non-smooth optimization one, whose
cost function is convex but not strictly or strongly convex.
In addition, the upper level optimization problem may appear indefinite
for some rather natural resource share allocations, which makes its solution
quite complicated.

In this paper, we propose to apply the exact non-smooth penalty method
together with the right-hand side decomposition and obtain a two-level
optimization problem, which yields a solution of the initial problem
under some fixed penalty parameter, unlike the smooth penalty method,
 and provides the consistency of constraints of
lower level problems for all the share allocations.
Nevertheless, for a successful implementation of this
decomposition method we have to solve the
upper level non-smooth convex optimization problem with a suitable simple
method since computation of any subgradient of the cost
function requires a solution of independent partial optimization problems.
We suggest a modification of the simplest subgradient projection method, which consists in
utilization of a two-speed step-size procedure. It does not require any a priori
information and provides the same minimal memory and computational expenses per iteration,
but enables us to generate step-sizes adaptively and reduce the number of improper steps.
We carried out series of computational experiments that confirmed
efficiency of the proposed methods.


\section{Decomposition via Shares Allocation} \label{s2}

We will utilize the following partition of the basic
$n$-dimensional Euclidean space
$$
\mathbb{R}^{n}=\mathbb{R}^{n_{1}}\times\mathbb{R}^{n_{2}}\times\ldots\times\mathbb{R}^{n_{l}},
$$
so that $n=\sum \limits_{i \in I} n_{i}$ with $I = \{1, \ldots , l\}$.
Similarly,  for each $x \in  \mathbb{R}^{n}$ we determine $x=(x_{i})_{i \in I}$
where $x_{i} = (x_{i1}, \dots,x_{i,n_{i}})^{\top}\in  \mathbb{R}^{n_{i}}$ for $i \in I$.
 Next, for each $i \in I$ let $X_{i}$ be a set in
 $\mathbb{R}^{n_{i}}$, $h_{i}: X_{i} \rightarrow
\mathbb{R}^{m}$ a mapping with components
$h_{ij} : X_{i} \rightarrow \mathbb{R}$ for $j=1, \ldots ,m$,
and $f_{i}: X_{i} \rightarrow \mathbb{R}$, $i \in I$ be
some functions.
We define the feasible set as follows:
\begin{equation} \label{eq:2.1}
 D=\left\{x \in X \ \vrule \ \sum _{i=1}^{l} h_{i}(x_{i}) \leq b \right\},
\end{equation}
where $X=X_{1} \times \dots \times X_{l}$,  $b $ is a
fixed vector in $\mathbb{R}^{m}$,
besides, let
\begin{equation} \label{eq:2.2}
 f (x) =  \sum _{i=1}^{l} f_{i} (x_{i}).
\end{equation}
We will take the following set of basic assumptions.

{\bf (A1)} {\em Each set $X_{i}$ is convex and closed and each
function $f_{i}: X_{i} \rightarrow \mathbb{R}$ is convex and continuous for $i\in I$.
Also, $h_{ij} : X_{i} \to \mathbb{R}$, $j=1,\ldots,m$,
$i\in I$ are convex continuous functions.}

Then (\ref{eq:1.1}), (\ref{eq:2.1}), and (\ref{eq:2.2})
give a general decomposable convex optimization problem.

Let us now introduce the set of  partitions
 of the right-hand side constraint vector $b$:
$$
U= \left\{u \in \mathbb{R}^{ml} \ \vrule \ \sum
\limits_{i=1}^{l}u_{i}  = b \right\},
$$
where $u= (u_{1},\ldots,u_{l})^{\top}$, $u_{i} \in \mathbb{R}^{m}$,
$i=1,\ldots,l$. Here  $u_{i}$ determines the $i$-th share of $b$.
The right-hand side decomposition method
is based on inserting an additional upper control level for finding the optimal shares
whereas the lower level optimization problem is decomposed into $l$
independent partial problems in $\mathbb{R}^{n_{i}}$ for $i \in I$; see \cite{KL65,Las70}.
More precisely, the approach consists first in the simple equivalent transformation of the feasible
set by inserting auxiliary variables:
$$
  D=\left\{x \in X \ \vrule \ \exists u \in U, \ h_{i}(x_{i})  \leq  u_{i}, \ i \in I\right \},
$$
Then problem (\ref{eq:1.1}), (\ref{eq:2.1}), and (\ref{eq:2.2}) is replaced by the
upper level optimization problem
\begin{equation} \label{eq:2.3}
 \min \limits _{u \in U} \to \tilde \mu(u)=  \sum _{i=1}^{l} \tilde \mu_{i} (u_{i}),
\end{equation}
where
\begin{equation} \label{eq:2.4}
\tilde \mu_{i} (u_{i})=\inf \{ f_{i} (x_{i}) \ | \ x_{i} \in X_{i}, \ h_{i}(x_{i})  \leq  u_{i}  \}
\end{equation}
is the marginal value function for the $i$-th partial problem, $i \in I$.

Observe that the initial optimization problem has $n$
variables and $m$ complex constraints, whereas problem
(\ref{eq:2.3}) has $ml$ variables and $m$
simple constraints, besides, computation of the value of each
function $\tilde \mu_{i} (u_{i})$ requires a solution of the optimization problem
having $n_{i}$ variables and $m$ simplified constraints. The preferences of
this substitution are obvious if $n$ is too large.
The main drawback of this approach is that the partial optimization problems in
(\ref{eq:2.4}) may have no solutions for some feasible partitions,
e.g. due to inconsistent constraints. Therefore, this approach
needs certain modifications.
For instance, it can be improved by applying a suitable penalty method.
Let $P_{i}(w_{i})=P_{i}(x_{i},u_{i})$ be a general
penalty function for the set
$$
W_{i}=\left\{w_{i}=(x_{i},u_{i}) \in \mathbb{R}^{n_{i}}\times \mathbb{R}^{m}\ \vrule \  h_{i}(x_{i})  \leq  u_{i} \right\},
$$
i.e.
$$
P_{i}(x_{i},u_{i}) \left \{
\begin{array}{cl}
=0, & \quad \mbox{if} \ (x_{i},u_{i}) \in W_{i}, \\
>0, &\quad \mbox{if} \ (x_{i},u_{i}) \notin W_{i};
\end{array}  \right.
$$
for $i \in I$. Then problem (\ref{eq:2.3})--(\ref{eq:2.4})
can be approximated by a sequence of penalized problems of the form
\begin{equation} \label{eq:2.5}
 \min \limits _{u \in U} \to \tilde \eta(u,\tau)=  \sum _{i=1}^{l} \tilde \eta_{i} (u_{i},\tau),
\end{equation}
where
\begin{equation} \label{eq:2.6}
\tilde \eta_{i} (u_{i},\tau)=\inf \{ f_{i} (x_{i}) +\tau P_{i}(x_{i},u_{i})\ | \ x_{i} \in X_{i} \}, \ i \in I,
\end{equation}
$\tau >0$ is a penalty parameter. Unlike (\ref{eq:2.4}), the partial optimization problems in
(\ref{eq:2.6}) always have solutions under general assumptions.
The above approach for separable convex optimization problems
was suggested in \cite{Raz67,Raz75,Umn75} where the known
smooth penalty functions were taken. Then the functions $\tilde \eta_{i}$
appear differentiable, however, finding a solution of problem
(\ref{eq:1.1}), (\ref{eq:2.1}), and (\ref{eq:2.2}) requires tending the parameter $ \tau$
to infinity, but large values of $ \tau$ give rather
difficult auxiliary problem (\ref{eq:2.5})--(\ref{eq:2.6}).
For this reason, we intend to apply the exact (non-smooth) penalty functions
$P_{i}(x_{i},u_{i})$ since penalized problem (\ref{eq:2.5})--(\ref{eq:2.6})
for some fixed $ \tau$ large enough then becomes equivalent to the initial problem.
This property of non-smooth penalties is well known for general optimization problems;
see \cite{Ere67}.


\section{Saddle Point Problems} \label{s3}

The suggested approach will be based on relationships with saddle point problems.
For this reason, we establish some properties of these problems associated with
the optimization problem (\ref{eq:1.1}), (\ref{eq:2.1}), and (\ref{eq:2.2}).
We define first its Lagrange function
$$
  M(x, \lambda)=\sum _{i=1}^{l} \left[f_{i}(x_{i})+\langle \lambda, h_{i}(x_{i})\rangle\right] -\langle \lambda, b\rangle,
$$
and the saddle point problem: Find a pair $(x^*,\lambda^*) \in X\times
\mathbb{R}^m_+$ such that
\begin{equation} \label{eq:3.1}
 \forall \lambda \in \mathbb{R}^m_+, \quad  M(x^*, \lambda)\leq M(x^*, \lambda^*)\leq M(x,
\lambda^*) \quad \forall x\in X,
\end{equation}
where
$$
\mathbb{R}^{m}_{+}=\left\{ {v \in \mathbb{R}^{m} \ | \ v_{i} \geq 0
 \quad i=1, \dots, m} \right\}.
$$
Problem (\ref{eq:3.1}) is equivalent to  the system:
\begin{eqnarray}
&& \displaystyle
\sum _{i=1}^{l} \left[f_{i}(x_{i})-f_{i}(x^*_{i})+\langle \lambda^*, h_{i}(x_{i})-h_{i}(x^*_{i})\rangle\right]\geq 0 \quad \forall x \in X,  \label{eq:3.2}\\
&& \displaystyle  \lambda^* \geq \mathbf{0}, \ b-\sum \limits_{i=1}^{l} h_{i}(x^*_{i}) \geq \mathbf{0}, \
\langle \lambda^*, b-\sum \limits_{i=1}^{l} h_{i}(x^*_{i})\rangle= 0.  \label{eq:3.3}
\end{eqnarray}
By using the suitable Karush-Kuhn-Tucker saddle point theorem for problem
(\ref{eq:1.1}), (\ref{eq:2.1}), and (\ref{eq:2.2})
(see e.g. \cite[Section 28]{Roc70}), we obtain that it is equivalent to
(\ref{eq:3.1}) under certain regularity type conditions.
So, together with {\bf (A1)} we will take the following basic assumption.

{\bf (A2)} {\em There exists a saddle point $(x^*,\lambda^*) \in X\times
\mathbb{R}^m_+$ in (\ref{eq:3.1}).}

Under these assumptions we have the minimax equality \begin{equation} \label{eq:3.4}
f^*=\inf_{x \in X} \sup_{\lambda \in \mathbb{R}^m_+}
  M(x, \lambda)=\sup_{\lambda \in \mathbb{R}^m_+} \inf_{x \in X}
  M(x, \lambda);
\end{equation}
see e.g. \cite[Ch. 1, Corollary 4.1]{GT89}.

Similarly, we can define the Lagrange function for
the same optimization problem with the auxiliary share variables
$$
  L(x,u,y)=\sum _{i=1}^{l} \left[f_{i}(x_{i})+\langle y_{i}, h_{i}(x_{i})-u_{i}\rangle\right],
$$
and the corresponding saddle point problem: Find a triple $(x^*,u^*,y^*) \in X\times U\times
\mathbb{R}^{ml}_{+}$ such that
\begin{equation} \label{eq:3.5}
\forall y \in \mathbb{R}^{ml}_{+}, \quad
L(x^*, u^*, y)\leq L(x^*, u^*, y^*)\leq L(x, u, y^*) \quad \forall x\in X, \ \forall u \in U.
\end{equation}
Problem (\ref{eq:3.5}) is equivalent to  the system:
\begin{eqnarray}
&& \displaystyle
\sum _{i=1}^{l} \left[f_{i}(x_{i})-f_{i}(x^*_{i})+\langle y_{i}^*, h_{i}(x_{i})-h_{i}(x^*_{i})+u^*_{i}-u_{i}\rangle\right]\geq 0 \quad \forall x \in X,  \label{eq:3.6}\\
&& \displaystyle  y_{i}^*= (1/l) \sum \limits_{s=1}^{l} y^*_{s}, \ i \in I,  \label{eq:3.7} \\
&& \displaystyle  y_{i}^* \geq \mathbf{0}, \ u^*_{i}- h_{i}(x^*_{i}) \geq \mathbf{0}, \
\langle y_{i}^*, u^*_{i}- h_{i}(x^*_{i})\rangle= 0, \ i \in I.  \label{eq:3.8}
\end{eqnarray}
Observe that (\ref{eq:3.7}) implies
\begin{equation} \label{eq:3.9}
y_{i}^*=y_{j}^* \quad \forall i\neq j,
\end{equation}
i.e. all the dual variables in (\ref{eq:3.5}) coincide.


\begin{proposition} \label{pro:3.1}
Suppose  {\bf (A1)} and {\bf (A2)} are fulfilled.

(i) If  a pair $(x^*,\lambda^*) \in X\times
\mathbb{R}^m_+$  is a saddle point in (\ref{eq:3.1}), then there exists a point $u^{*} \in U$ such that
the triple $(x^*,u^*,y^*)$ where $y_{i}^*= \lambda^*$, $i \in I$, is a saddle point in (\ref{eq:3.5}).

(ii) If a triple $(x^*,u^*,y^*)  \in X\times U\times
\mathbb{R}^{ml}_{+}$ is a saddle point in (\ref{eq:3.5}), then
a pair $(x^*,\lambda^*)$  is a saddle point in (\ref{eq:3.1}) where $\lambda^*=y_{i}^*$ for any $i \in I$.
\end{proposition}
{\bf Proof.}
Let $(x^*,\lambda^*) \in X\times
\mathbb{R}^m_+$  be a saddle point in (\ref{eq:3.1}). Set
$\tilde u_{i}=h_{i}(x^*_{i})$, $i \in I$. Then (\ref{eq:3.3}) gives
$\sum \limits_{i=1}^{l} \tilde u_{i} \leq b$. Now define
$$
u^*_{i}= (1/l) \left[b-\sum \limits_{s=1}^{l} \tilde u_{s} \right]\ +\tilde u_{i}
 \ \mbox{and} \ y_{i}^*= \lambda^*, \ i \in I,
$$
then (\ref{eq:3.7}) holds and (\ref{eq:3.3}) gives (\ref{eq:3.8}).
Also, (\ref{eq:3.2}) gives (\ref{eq:3.6}). Hence,
$(x^*,u^*,y^*)$ is a saddle point in (\ref{eq:3.5}).

Conversely, let $(x^*,u^*,y^*)$ be a saddle point in (\ref{eq:3.5}).
Due to (\ref{eq:3.9}), we can set $\lambda^*=y_{i}^*$ for any $i \in I$.
Also, (\ref{eq:3.6}) gives (\ref{eq:3.2}), whereas (\ref{eq:3.8}) gives (\ref{eq:3.3}).
Hence, $(x^*,\lambda^*)$ is a saddle point in (\ref{eq:3.1}).
\QED

Combining the above assertions and (\ref{eq:3.4}) we also obtain the minimax equality \begin{equation} \label{eq:3.10}
f^*=\inf_{x \in X, u \in U} \sup_{y \in \mathbb{R}^{ml}_{+}}
  L(x,u,y)=\sup_{y \in \mathbb{R}^{ml}_{+}} \inf_{x \in X,u \in U}
  L(x,u,y).
\end{equation}

We can take a parametric vector $t \in \mathbb{R}^{m}_{+}$ and define the reduced set of
dual variables:
$$
Y_{t}= \left\{y \in \mathbb{R}^{ml}_{+} \ \vrule \
y_{i} \leq t, \ i \in I \right\}.
$$
Then, by analogy with (\ref{eq:3.5}),  we can define the
 modified saddle point problem: Find a triple $(x^*,u^*,y^*) \in X\times U\times
Y_{t}$ such that
\begin{equation} \label{eq:3.11}
\forall y \in Y_{t}, \quad
L(x^*, u^*, y)\leq L(x^*, u^*, y^*)\leq L(x, u, y^*) \quad \forall x\in X, \ \forall u \in U.
\end{equation}
If $(x^*,\lambda^*) \in X\times
\mathbb{R}^m_+$  is a saddle point in (\ref{eq:3.1}) and $\lambda^* \leq t$, then it follows from
Proposition \ref{pro:3.1} that the corresponding
triple $(x^*,u^*,y^*)$, which is a saddle point in (\ref{eq:3.5}), where $y_{i}^*= \lambda^*$, $i \in I$,
will be a saddle point in (\ref{eq:3.11}). Moreover, (\ref{eq:3.10}) now implies
 the minimax equality \begin{equation} \label{eq:3.12}
f^*=\inf_{x \in X, u \in U} \sup_{y \in Y_{t}}
  L(x,u,y)=\sup_{y \in Y_{t}} \inf_{x \in X,u \in U}
  L(x,u,y).
\end{equation}
Therefore, saddle points in (\ref{eq:3.11})  are also related to solutions of
the optimization problem (\ref{eq:1.1}), (\ref{eq:2.1}), and (\ref{eq:2.2}).


\begin{proposition} \label{pro:3.2}
Let {\bf (A1)} and {\bf (A2)} be fulfilled and let $(x^*,\lambda^*) \in X\times
\mathbb{R}^m_+$  be a saddle point in (\ref{eq:3.1}).
If  $\lambda^* \leq t$, then there exists a point $u^{*} \in U$ such that
the triple $(x^*,u^*,y^*)$ where $y_{i}^*= \lambda^*$, $i \in I$, is a saddle point in
(\ref{eq:3.5}) and (\ref{eq:3.11})
and the minimax equality (\ref{eq:3.12}) holds true.
\end{proposition}


\section{Exact Decomposable Penalty Method} \label{s4}

Let us select the primal optimization problem in (\ref{eq:3.12}):
\begin{equation} \label{eq:4.1}
\min \limits _{x \in X, u \in U} \to  \sup_{y \in Y_{t}}
  L(x,u,y).
\end{equation}
Since
\begin{eqnarray*}
 \displaystyle
\sup_{y \in Y_{t}}
  L(x,u,y) &=& \sup_{y \in Y_{t}} \sum _{i=1}^{l} \left[f_{i}(x_{i})+\langle y_{i}, h_{i}(x_{i})-u_{i}\rangle\right] \\
 \displaystyle  &=&  \sum _{i=1}^{l} \max_{\mathbf{0} \leq y_{i}  \leq t} \left[f_{i}(x_{i})+\langle y_{i}, h_{i}(x_{i})-u_{i}\rangle\right] \\
 \displaystyle  &=&  \sum _{i=1}^{l} \left[f_{i}(x_{i})+P_{i}(x_{i},u_{i},t)\right],
\end{eqnarray*}
where
$P_{i}(x_{i},u_{i},t)=\langle t, [h_{i}(x_{i})-u_{i}]_{+}\rangle$ for $i \in I$,
$[v]_{+}$ denotes the projection of $v \in \mathbb{R}^{m}$ onto the non-negative orthant
$\mathbb{R}^{m}_{+}$, we can rewrite (\ref{eq:4.1}) as follows:
\begin{equation} \label{eq:4.2}
\min \limits _{x \in X, u \in U} \to \sum _{i=1}^{l} \left[f_{i}(x_{i})+P_{i}(x_{i},u_{i},t)\right].
\end{equation}
This is nothing but the non-smooth penalty method problem  for the initial
optimization problem (\ref{eq:1.1}), (\ref{eq:2.1}), and (\ref{eq:2.2})
with the auxiliary share variables. We now give the basic equivalence properties.


\begin{theorem} \label{thm:4.1}
Let {\bf (A1)} and {\bf (A2)} be fulfilled and let $(x^*,\lambda^*) \in X\times
\mathbb{R}^m_+$  be a saddle point in (\ref{eq:3.1}).

(i) If  $\lambda^* \leq t$, then $x^{*}$  is a solution of
problem (\ref{eq:4.2}), which has the optimal value $f^{*}$.

(ii) If  $\lambda^* < t$, then any solution $ (\bar x,\bar u)$  of
problem (\ref{eq:4.2}) solves also
problem (\ref{eq:1.1}), (\ref{eq:2.1}), and (\ref{eq:2.2}).
\end{theorem}
{\bf Proof.}
Part (i) follows directly from Proposition \ref{pro:3.2}.
Next, set
$$
  W=\left\{(x,u) \in X\times U \ \vrule \  h_{i}(x_{i})  \leq  u_{i}, \ i \in I\right \}
$$
and take any $t \in \mathbb{R}^m_+$ so that problem (\ref{eq:4.2}) has a solution pair $(x(t),u(t))$. Then
\begin{eqnarray}
 \displaystyle
  f(x(t)) &\leq & \sum _{i=1}^{l} \left[f_{i}(x_{i}(t))+P_{i}(x_{i}(t),u_{i}(t),t)\right] \nonumber \\
 \displaystyle  &\leq & \inf \limits _{(x,u) \in W}  \sum _{i=1}^{l} \left[f_{i}(x_{i})+P_{i}(x_{i},u_{i},t)\right] =  \inf \limits _{x \in D}  f(x)=f^{*}. \label{eq:4.3}
\end{eqnarray}
Take now any $t' > \lambda^*$, then there exists $t''$ such that $t' > t'' > \lambda^*$.
For brevity, set $x' =x(t')$, $u' =u(t')$ and $x'' =x(t'')$, $u'' =u(t'')$, these elements
exist due to (i). From Proposition \ref{pro:3.2} it follows that there exists a point $u^{*} \in U$ such that
the triple $(x^*,u^*,y^*)$ where $y_{i}^*= \lambda^*$, $i \in I$, is a saddle point in
(\ref{eq:3.5}) and (\ref{eq:3.11}) both for $t=t'$ and $t=t''$, moreover,
the minimax equalities (\ref{eq:3.10}) and (\ref{eq:3.12}) both for $t=t'$ and $t=t''$ hold true.
Combining these relations together with (\ref{eq:4.3}) we obtain
\begin{eqnarray*}
 \displaystyle
  && \sum _{i=1}^{l} \left[f_{i}(x'_{i})+P_{i}(x'_{i},u'_{i},t')\right] \leq f^{*} = L(x^*, u^*, y^*) \\
 \displaystyle  && \leq L(x', u', y^*) \leq \sup_{y \in Y_{t''}} L(x', u', y)
 = \sum _{i=1}^{l} \left[f_{i}(x'_{i})+P_{i}(x'_{i},u'_{i},t'')\right].
\end{eqnarray*}
It follows that
$$
\sum _{i=1}^{l} \langle t'-t'', [h_{i}(x'_{i})-u'_{i}]_{+}\rangle \leq 0,
$$
hence, $(x',u') \in W$ and $x' \in D$. Due to (\ref{eq:4.3}), this gives $x' \in D^{*}$.
\QED

These results can be viewed as a specialization the known properties of non-smooth penalty functions
in general optimization problems (see \cite{Ere67}) to decomposable optimization problems.

Since the variables $x$ and $u$ in (\ref{eq:4.2}) are not contained in joint constraints,
we can apply the sequential minimization and take
the equivalent optimization problem
\begin{equation} \label{eq:4.4}
 \min \limits _{u \in U} \to \mu(u,t)=  \sum _{ i \in I} \mu_{i} (u_{i},t),
\end{equation}
where
\begin{equation} \label{eq:4.5}
\mu_{i} (u_{i},t)=\inf \limits _{x_{i} \in X_{i}} \{ f_{i}(x_{i})+P_{i}(x_{i},u_{i},t)\}, \ i \in I,
\end{equation}
which gives precisely the exact (non-smooth) decomposable penalty function method; cf.
(\ref{eq:2.3})--(\ref{eq:2.4}) and (\ref{eq:2.5})--(\ref{eq:2.6}).
Now Theorem \ref{thm:4.1} guarantees that a solution of (\ref{eq:4.4})--(\ref{eq:4.5})
yields a solution of (\ref{eq:2.3})--(\ref{eq:2.4}) if $\lambda^* < t$, but now the $i$-th partial problem
(\ref{eq:4.5}) does not contain inconsistent constraints for any partition
of the right-hand side vector $b$. As in (\ref{eq:2.5})--(\ref{eq:2.6}),
calculation of the marginal value function for each partial problem
can be made independently and separately of the others. Hence,
we have derived another two-level decomposition method
that has certain preferences over the
other right-hand side decomposition methods.

Under the assumptions in  {\bf (A1)} the functions $\mu_{i} $, $i \in I$, are convex in $u_{i}$
(see e.g. \cite[Ch. 1, Theorem 5.11]{GT89}), hence so is $\mu$ in $u$. At the same time,
$\mu$ need not be differentiable in $u$ in general. Its subdifferential
$\partial_{u} \mu (u,t)$ in $u$ can be found by the proper formula for the composite
convex functions; see \cite[Ch. 3, Theorem 2.6]{Psh80} and \cite[Ch. 1, Theorem 5.11]{GT89}.
However, we now give a specialization of this formula for calculations of a subgradient for
each marginal value function in (\ref{eq:4.5}), which is more suitable for
utilization in iterative solution methods for problem (\ref{eq:4.4})--(\ref{eq:4.5}).
In fact, due to the minimax equality we have
\begin{eqnarray*}
 \displaystyle
   \mu_{i} (u_{i},t)&=&\inf \limits _{x_{i} \in X_{i}} \{ f_{i}(x_{i})+P_{i}(x_{i},u_{i},t)\}
   = \inf \limits _{x_{i} \in X_{i}} \sup_{\mathbf{0} \leq y_{i}  \leq t} \{f_{i}(x_{i})+\langle y_{i},h_{i}(x_{i})-u_{i}\rangle\} \\
 \displaystyle  &=& \sup_{\mathbf{0} \leq y_{i}  \leq t} \inf \limits _{x_{i} \in X_{i}} \{f_{i}(x_{i})+\langle y_{i},h_{i}(x_{i})-u_{i}\rangle\} \\
 \displaystyle  &=& \sup_{\mathbf{0} \leq y_{i}  \leq t} \left\{ \inf \limits _{x_{i} \in X_{i}}
\left[ f_{i}(x_{i})+\langle y_{i},h_{i}(x_{i})\rangle\right] - \langle y_{i},u_{i}\rangle\right\}
= \sup_{\mathbf{0} \leq y_{i}  \leq t} \sigma_{i}(y_{i},u_{i}).
\end{eqnarray*}
Take arbitrary points $u'_{i}, u''_{i} \in \mathbb{R}^{m}$. Let
$$
y'_{i}=\arg \max \limits_{\mathbf{0} \leq y_{i}  \leq t}\sigma_{i}(y_{i},u'_{i}),
y''_{i}=\arg \max \limits_{\mathbf{0} \leq y_{i}  \leq t}\sigma_{i}(y_{i},u''_{i}).
$$
Then
\begin{eqnarray*}
 \displaystyle
   \mu_{i} (u''_{i},t)-\mu_{i} (u'_{i},t)&=& \sigma_{i}(y''_{i},u''_{i})-\sigma_{i}(y'_{i},u'_{i})
   \geq \sigma_{i}(y'_{i},u''_{i})-\sigma_{i}(y'_{i},u'_{i}) \\
   &=& \langle -y'_{i},u''_{i}-u'_{i}\rangle.
\end{eqnarray*}
This means that $-y'_{i} \in \partial_{u_{i}} \mu_{i} (u'_{i},t)$.
Therefore, the calculation of some subgradient of $\mu_{i}(u'_{i},t)$ in $u_{i}$
reduces to finding a solution of the problem
\begin{equation} \label{eq:4.6}
 \max  \limits_{\mathbf{0} \leq y_{i}  \leq t} \to \left\{ \inf \limits _{x_{i} \in X_{i}}
\left[ f_{i}(x_{i})+\langle y_{i},h_{i}(x_{i})\rangle\right] - \langle y_{i},u'_{i}\rangle\right\},
\end{equation}
which is precisely the modified dual optimization problem to (\ref{eq:2.4}) with the additional upper bound
$y_{i}  \leq t$.


\section{Step-size Strategies for Subgradient Projection Methods} \label{s5}

Together with the above example there exist a great number of some other significant
applications of convex minimization problems containing just
non-differentiable functions; see \cite{Las70,Erm76,Sho79,Pol83,Pan85,MN92} and the references therein.
For this reason, their theory and methods were developed rather well.
We recall that most applications admit calculation of only one arbitrary taken element from the
subdifferential of a non-differentiable function at any point.
During a rather long time, most efforts
were concentrated on developing more powerful and rapidly convergent methods
such as space dilation and bundle type ones.
However, significant areas of applications related to decision making in
industrial, transportation, information and communication systems,
having large dimensionality and inexact data together with
scattered necessary information force one to avoid
complex transformations and even line-search procedures, which are
involved in all the mentioned methods. Let us take
the problem of minimizing a convex, but not necessarily differentiable
function  $\varphi : E \to \mathbb{R}$
 on a convex set $V \subseteq E$ in a finite-dimensional
 Euclidean space $E$, or briefly,
\begin{equation} \label{eq:5.1}
 \min \limits _{v\in V} \to \varphi(v).
\end{equation}
Its solution set is denoted by $V^{*}$ and the optimal value of the
function by $\varphi^{*}$. Next, the constraint set is supposed to be rather simple
in the sense that the projection of a point $x$ onto $V$, which is denoted
by $\pi _{V} (x)$, is not very expensive.

Then one can apply the simplest and most popular subgradient projection method,
which provides convergence to a solution under the so-called divergent series step-size rule.
Its iteration computation expenses
and accuracy requirements are rather low, but its convergence may be rather slow;
see e.g. \cite[Ch. 2, \S \S  1--2]{Sho79} and
\cite[Ch. 5, \S 3]{Pol83}.
There are several ways to speed up convergence of the subgradient methods
via utilization of a priori information such as the optimal value
or some condition numbers.
However, it is usually difficult to calculate these values exactly,
whereas taking inexact estimates may again lead to slow convergence.
For this reason, we will try to improve convergence of
the subgradient projection method within the divergent series step-size rule.

In this section, we will take the following assumptions.

{\bf (B1)} {\em The set $V \subseteq E$ is convex and closed and the
function $\varphi: V \rightarrow \mathbb{R}$ is convex and continuous.
Also, the set $V^{*}$ is nonempty. }

{\bf (B2)} {\em There exist a number $C< +\infty$ such that}
$$
\|g\| \leq C, \quad \forall g\in \partial  \varphi (v), \ \forall v\in V.
$$

Now we write the simplest subgradient projection method for problem (\ref{eq:5.1}):
\begin{equation} \label{eq:5.2}
 v^{k+1}=\pi _{V} [v^k-\theta _{k} g^k], \quad g^k \in
\partial\varphi(v^{k}),
\end{equation}
where
\begin{equation} \label{eq:5.3}
\theta _{k} > 0, \quad \sum \limits_{k=0}^{\infty } \theta _{k} = \infty , \quad
  \sum \limits_{k=0}^{\infty } \theta_{k}^{2} < \infty.
\end{equation}
The method stops if $g^k=\mathbf{0}$ or $v^{k+1}=v^k$,
then $v^k \in V^{*}$, but we suppose this situation does not occur.
Its convergence properties can be formulated as follows
(see \cite[Ch. 1, \S 3, p.47]{Erm76}, \cite[Ch. 3, Theorem 4.5]{DV81}, and
\cite[Ch. 5, \S 2]{GT89}).


\begin{proposition} \label{pro:5.1}
Let {\bf (B1)} and {\bf (B2)} be fulfilled and let
the sequence $\{ v^{k} \}$ be generated in conformity with (\ref{eq:5.2})--(\ref{eq:5.3}).
Then
$$
 \lim _{k \rightarrow \infty } v^{k} = v^{*} \in V^{*}.
$$
\end{proposition}
Besides, if we replace the subgradient $g^k$ in (\ref{eq:5.2})
with the normed subgradient $q^k=(1/\|g^k\|)g^k$, the assertion of
Proposition \ref{pro:5.1} remains true without {\bf (B2)}.

The slow convergence of the sequences generated by method
 (\ref{eq:5.2})--(\ref{eq:5.3}) is mainly due to
the divergent series condition in (\ref{eq:5.3}) that is
non-adaptive to iterates since nonlinear functions may behave in a different manner on
different parts of their domains. Usually, the sequence of step-sizes $\{\theta _{k} \}$
has a unique rate of decrease, e.g.
$$
  \theta _{k}=\theta /(k+1)^{\tau}, \ \theta > 0,
$$
where $\tau \in (0.5,1]$. Then the sequence $\{\theta _{k} \}$ may contain many
rather large improperly chosen step-sizes that do not in fact decrease
the distance to a solution and cause the slow convergence.
In order to reduce the number of these improper step-sizes
we propose to apply a two-speed step-size procedure within (\ref{eq:5.3}).

Namely, choose an index sequence $\{i _{s}\}$ such that
\begin{equation} \label{eq:5.5}
 i _{0}=0, \ 0 < i _{s+1}-i _{s} \leq d< \infty, \  s=0,1, \ldots ,
\end{equation}
and a sequence $\{\beta _{s}\}$ such that
\begin{equation} \label{eq:5.6}
\beta _{s} > 0, \quad \sum \limits_{s=0}^{\infty } \beta _{s} = \infty , \quad
  \sum \limits_{s=0}^{\infty } \beta _{s}^{2} < \infty.
\end{equation}
Then apply the subgradient projection method (\ref{eq:5.2}) where
\begin{equation} \label{eq:5.7}
\theta _{k}=\left \{
\begin{array}{cl}
\beta _{s}, & \quad \mbox{if} \ k=i _{s}, \\
\nu\theta _{k-1}, &\quad \mbox{if} \ i _{s}<k< i _{s+1},
\end{array}  \right. \quad \forall s=0,1, \ldots, \ \nu \in (0,1).
\end{equation}

It follows that the more rapid step-size decrease at iterates between $i _{s}$
 and $i _{s+1}$ allows us to find more suitable steps, but
 the divergent series condition in (\ref{eq:5.6}) prevents from the very small steps.
 It should be noted that the proposed two-speed step-size procedure
 differs from the relaxation type versions of the subgradient projection method
 (see \cite[Ch. 3, \S 4]{DV81}) since our version does not provide any
monotone decrease of the goal function at each iteration and does not involve any line-search.

At the same time, it is not difficult to show that  (\ref{eq:5.5})--(\ref{eq:5.7})
imply (\ref{eq:5.3}). In fact,
$$
\sum \limits_{k=0}^{\infty } \theta _{k} \geq \sum \limits_{s=0}^{\infty } \beta _{s} = \infty
$$
and
$$
\sum \limits_{k=0}^{\infty } \theta_{k}^{2} \leq d \sum \limits_{s=0}^{\infty } \beta _{s}^{2} < \infty .
$$
Therefore, the assertion of Proposition \ref{pro:5.1} holds true for
the subgradient projection method (\ref{eq:5.2}), (\ref{eq:5.5})--(\ref{eq:5.7}).

Next, if we write problem (\ref{eq:4.4})--(\ref{eq:4.5})  in the format (\ref{eq:5.1}), then
{\bf (A1)}--{\bf (A2)} imply {\bf (B1)}--{\bf (B2)}. Since the partial subgradients are calculated
as negative solutions of problems (\ref{eq:4.6}), we can take $C=\sqrt{l}\|t\|$.


\section{Computational Experiments}\label{s6}

In order to check the performance of the proposed methods we carried
out preliminary series of computational experiments.
We chose two classes of test problems and implemented all the methods
mentioned below in Delphi with double precision
arithmetic.


\subsection{Non-smooth test problem}\label{s6.1}

Our first goal was to compare method (\ref{eq:5.2}), (\ref{eq:5.5})--(\ref{eq:5.7})
with other simple subgradient methods.
We took for comparison the well-known non-smooth
test problem of form (\ref{eq:5.1}) with $V=E$ from \cite{ShS72}, where
$$
\varphi(v) = \max \limits _{i=1, \dots ,m} \eta_{i}(v), \
 \eta_{i}(v) = b_{i} \sum \limits _{j=1}^{n}(v_{j} - a_{ij})^{2}, \ i=1, \dots ,m,
$$
with $n=5$ and $m=10$. The coefficients of the quadratic functions are given by the vector
$b = (1, 5, 10, 2, 4, 3, 1.7, 2.5, 6, 3.5)^{\top}$ and by the transposed matrix
$$
A^{\top}=  \left( \begin{array}{ccccccccccc}
0 & 2 & 1 & 1 & 3 & 0 & 1 & 1 & 0 & 1 \\
0 & 1 & 2 & 4 & 2 & 2 & 1 & 0 & 0 & 1 \\
0 & 1 & 1 & 1 & 1 & 1 & 1 & 1 & 2 & 2 \\
0 & 1 & 1 & 2 & 0 & 0 & 1 & 2 & 1 & 0 \\
0 & 3 & 2 & 2 & 1 & 1 & 1 & 1 & 0 & 0
    \end{array} \right);
$$
see also \cite{Sho79,Pol83}. The problem has
the optimal value $\varphi^{*}=22.60016$. We chose
the starting point $v^{0} = (0,0,0,0,1)^{\top} $.

We chose the simplest subgradient method (\ref{eq:5.2}) where
$$
\theta _{k}=\theta/(k+1), \quad \theta>0.
$$
It was abbreviated as (SGM). For method (\ref{eq:5.2}), (\ref{eq:5.5})--(\ref{eq:5.7}),
which was abbreviated as (SGMTS), we chose the rules
$$
\beta _{s}=\theta/(s+1), \ \theta>0, \ i _{s+1}-i _{s} = d.
$$
We also took the same method (\ref{eq:5.2}) where
\begin{equation} \label{eq:6.3}
\theta _{k}=\theta/\sqrt{(k+1)}, \quad \theta>0.
\end{equation}
It was abbreviated as (SGMSQ).
In addition, the so-called
subgradient method of simple double averaging from \cite{NeS15} was used. It  can be written as follows:
\begin{eqnarray*}
 & &  v^{k+1} =\mu _{k}v^{k}+ (1 -\mu _{k}) y^{k}, \ \mu _{k}=(k+1)/(k+2),     \\
  & &   y^{k} =v^{0} - \theta _{k} p^{k}, \ p^{k}=\sum \limits_{i=0}^{k } g^{i}, \
 g^{i} \in \partial \varphi (v^{i}),
\end{eqnarray*}
where $ \theta _{k}$ was chosen as in (\ref{eq:6.3}).
We abbreviate this method as (DASG).

We compared all the methods for different accuracy $\varepsilon$
with respect to the goal function deviation
$$
\Delta (v)=\varphi(v)-\varphi^{*}.
$$
We took the same starting step $\theta =0.1$ for all the methods.
For (SGMTS), we took the ratio $\nu = 0.7$ and set $d = 25$.
The results  are given in Table \ref{tbl:1}, where we indicate
the total number of iterations (it)
(or the total number of subgradient calculations) for attaining
the desired accuracy $\varepsilon$.

\begin{table} \caption{Comparison of subgradient methods} \label{tbl:1}
\centering
\begin{tabular}[t]{|r|r|r|r|r|r|r|r|}
\hline
\multicolumn{2}{|c|}{SGM} & \multicolumn{2}{|c|}{SGMTS} & \multicolumn{2}{|c|}{SGMSQ} & \multicolumn{2}{|c|}{DASG}\\
\hline
   {$\varepsilon$} & {it} & {$\varepsilon$} & {it} & {$\varepsilon$} & {it} & {$\varepsilon$} & {it}\\					
\hline
0.1   & 60       & 0.1 & 21      & 0.1 & 404       & 0.1 & 117\\					
\hline
 0.01 & 252      & 0.01 & 292    & 0.01 & 14575    & 0.01 & 1542\\
\hline
 0.001 & 1410    & 0.001 & 570   & 0.003 & 35000  & 0.001 & 9982\\	
\hline
 0.0001 & 6728  & 0.0001   & 3696 & - & -          & 0.0003 & 35000\\	
\hline
\end{tabular}
\end{table}

The implementation of (SGMTS) showed rather rapid convergence
in comparison with the other methods.
Convergence of (SGM) appeared rather slow, but stable. At the same time,
(SGM) showed better convergence properties than the subgradient method
utilizing rule (\ref{eq:6.3}), even after the averaging procedure.


\subsection{Decomposable linear programming test problems}\label{s6.2}

Afterwards, we took decomposable linear programming
test examples of form (\ref{eq:1.1}), (\ref{eq:2.1}), and (\ref{eq:2.2}). The main goal
was to investigate convergence of the subgradient projection methods to a solution of
the non-smooth master problem (\ref{eq:4.4})--(\ref{eq:4.5}).
More precisely, the initial problem was the following:
$$
\max \ \to \ \left\{ \sum \limits_{i=1}^{l} \langle c_{i}, x_{i} \rangle \ \vrule \
\sum \limits_{i=1}^{l} A_{i}x_{i}  \leq  b, \ x_{i} \geq \mathbf{0}, \ i = 1, \ldots , l \right\}.
$$
It can be treated as a model describing the economic system which deals in $n$ commodities
and $m$ pure factors of production. These common factors are utilized by
$l$ producers, so that the $i$-th
producer chooses an output vector $x_{i}\in \mathbb{R}^{n_{i}}$,
his/her consumption rates are described by an $m \times n_{i}$ matrix
$A_{i}$, whereas the vector $c_{i}\in \mathbb{R}^{n_{i}}$ denotes prices of
his/her  outputs, the vector $b\in \mathbb{R}^{m}$ denotes inventories of common
factors. Therefore, the system should choose the outputs for maximizing the income value.

This problem is rewritten equivalently as follows:
$$
\min \ \to \ \left\{ \sum \limits_{i=1}^{l} \langle -c_{i}, x_{i} \rangle \ \vrule \
A_{i}x_{i}  \leq  u_{i}, \ x_{i} \geq \mathbf{0}, \ i = 1, \ldots , l, \ u \in U \right\},
$$
which corresponds to the format (\ref{eq:2.3})--(\ref{eq:2.4}), where
$$
f_{i} (x_{i})=\langle -c_{i}, x_{i} \rangle, \ h_{i}(x_{i})=A_{i}x_{i}, \ X_{i}=\mathbb{R}^{n_{i}}_{+}, \ i = 1, \ldots , l.
$$
Hence, we can  apply the exact non-smooth decomposable penalty method to this problem
and solve its basic problem (\ref{eq:4.4})--(\ref{eq:4.5}) with a suitable
subgradient projection method. In this case, the calculation of some subgradient of the function
$\mu_{i}(u'_{i},t)$ in $u_{i}$
reduces to finding a solution of the partial linear programming problem
$$
\min \ \to \ \left\{  \langle u_{i}, y_{i} \rangle \ \vrule \
A^{\top}_{i}y_{i}  \geq  c_{i}, \ \mathbf{0} \leq y_{i} \leq t \right\}.
$$
Next, iterate (\ref{eq:5.2}) with $V=U$ is
rewritten equivalently as follows:
\begin{equation} \label{eq:6.5}
 v^{k+1}=v^k-\theta _{k} \bar g^k, \quad \bar g^k_{i}= g^k_{i}- (1/l) \sum \limits_{i=1}^{l}
g^k_{s}, \  i = 1, \ldots , l, \quad g^k \in \partial\varphi(v^{k}).
\end{equation}
We chose $m=2$ and $n_{i}=2$
for each $i=1,\ldots ,l$ and took different values of $l$. Then the
master problem  (\ref{eq:2.5}) has $2l$ variables. The elements of the above matrices and vectors
were generated to be positive with the help of trigonometric functions. The upper bound $t$
was somewhat different for varying elements. In all the
cases, we took the same starting point $u^{0}$ with $u^{0}_{i}=(1/l)b$ for $i \in I$.

Since  methods (SGMSQ) and (DASG) appeared too slow, we chose the subgradient projection
method (\ref{eq:6.5}) with the rule
$$
\theta _{k}=\theta/(k+2), \quad \theta>0,
$$
which was also abbreviated as (SGM), and method (\ref{eq:6.5}), (\ref{eq:5.5})--(\ref{eq:5.7}),
which was abbreviated as (SGMTS), where we chose  the rules
$$
\beta _{s}=\theta/(s+2), \ \theta>0, \ i _{s+1}-i _{s} = d.
$$
We fixed the value $ \theta =5$ for both the methods. Since
the optimal values of test problems were unknown
we show the number of iterations (it) and the best attained value of the goal function (f).
The results  are given in Tables \ref{tbl:2} and \ref{tbl:3}.

\begin{table} \caption{Computations by (SGM)} \label{tbl:2}
\centering
\begin{tabular}[t]{|r|r|r|r|r|r|r|r|}
\hline
\multicolumn{2}{|c|}{$l=2$} & \multicolumn{2}{|c|}{$l=10$} & \multicolumn{2}{|c|}{$l=20$} & \multicolumn{2}{|c|}{$l=50$}\\
\hline
    {it} & {$f$} & {it} & {$f$} & {it} & {$f$} & {it} & {$f$} \\						
\hline
0   & -6.3536       & 0   & -6.5946     & 0   & -7.8872       & 0   & -7.9186\\					
\hline
50 & -7.3354      & 50 & -14.0226    & 50 & -12.8349     & 50 & -11.928 \\
\hline
 100 & -          & 100 & -14.6558   & 100 & -14.911  & 100 & -14.6894\\	
\hline
 200 & -           & 400   & -15.1555 & 350 & -16.2088 & 1950 & -16.7793\\	
\hline
\hline
{$t_{1}$} & {$t_{2}$} & {$t_{1}$} & {$t_{2}$} & {$t_{1}$} & {$t_{2}$} & {$t_{1}$} & {$t_{2}$} \\						
\hline
2.62   & 3.92       & 2.62   & 6.01     & 2.66   & 6.37       & 2.67   & 6.41  \\					
\hline
\end{tabular}
\end{table}

The implementation of (SGMTS) for different dimensionalities showed also more rapid convergence
in comparison with  (SGM).

\begin{table} \caption{Computations by (SGMTS)} \label{tbl:3}
\centering
\begin{tabular}[t]{|r|r|r|r|r|r|r|r|}
\hline
\multicolumn{2}{|c|}{$l=2$} & \multicolumn{2}{|c|}{$l=10$} & \multicolumn{2}{|c|}{$l=20$} & \multicolumn{2}{|c|}{$l=50$}\\
\hline
{$\nu = 0.2$} & {$d = 10$} & {$\nu = 0.8$} & {$d = 25$} & {$\nu = 0.9$} & {$d = 40$} & {$\nu = 0.9$} & {$d = 100$} \\
\hline
    {it} & {$f$} & {it} & {$f$} & {it} & {$f$} & {it} & {$f$} \\						
\hline
0   & -6.3536       & 0   & -6.5946     & 0   & -7.8872       & 0   & -7.9186\\					
\hline
50 & -7.0786      & 50 & -14.8631       & 50 & -14.7241       & 50 & -13.058 \\
\hline
 100 & -7.2158    & 100 & -15.952       & 100 & -15.2835       & 100 & -14.7465\\	
\hline
 200 & -7.3739    & -   & -             & 150 & -18.519        & 250 & -19.5964\\	
\hline
\hline
{$t_{1}$} & {$t_{2}$} & {$t_{1}$} & {$t_{2}$} & {$t_{1}$} & {$t_{2}$} & {$t_{1}$} & {$t_{2}$} \\						
\hline
2.62   & 3.92       & 2.62   & 6.01     & 2.66   & 6.37       & 2.67   & 6.41  \\					
\hline
\end{tabular}
\end{table}


\section*{Acknowledgement}

In this work, the author was supported by Russian Foundation for Basic Research, project No.
19-01-00431.



\begin{thebibliography}{99}
\markboth{I.V.~Konnov}{Exact penalties for optimization problems}

\bibitem
{Las70}
L.S.~Lasdon, {\em Optimization Theory for Large Systems}, Macmillan, New York, 1970.

\bibitem{Tsu81}
V.I.~Tsurkov,
{\em Decomposition in Large-Scale Problems},  Nauka, Moscow, 1981. [In Russian]

\bibitem{GT89}
E.G.~Gol'shtein and N.V.~Tret'yakov,
{\em Augmented Lagrange Functions}, Nauka,
Moscow, 1989 (Engl. transl. in John Wiley and Sons, New York, 1996).

\bibitem{CCM06}
A.J.~Conejo, E.~Castillo, R.~Minguez, and R.~Garcia-Bertrand,
{\em Decomposition Techniques in Mathematical Programming},
Springer-Verlag, Berlin, 2006.

\bibitem{KL65}
J.~Kornai and T.~Liptak, {\em Two-level planning}, Econometrica, vol.33
(1965), pp.141--169.

\bibitem{Raz67}
B.S.~Razumikhin, {\em Iterative method for the solution
and decomposition of linear programming problems},
 Autom. Remote. Control., vol.29 (1967), pp.427--443.

\bibitem{Raz75}
B.S.~Razumikhin,  {\em  Physical Models and Methods of Equilibrium
Theory in Programming and Economics}, Nauka, Moscow, 1975. [In Russian]

\bibitem{Umn75}
 A.E.~Umnov, {\em The method of penalty functions in problems of large
dimension},  USSR Comp. Maths. Math. Phys., vol.15 (1975), pp.32--45.

\bibitem{Ere67}
I.I.~Eremin, {\em The penalty method in convex programming},  Cybernetics, vol.3 (1967), pp. 53--56.

\bibitem{Roc70}
R.T.~Rockafellar, {\em Convex Analysis},
Princeton University Press, Princeton, 1970.

\bibitem{Psh80}
B.N.~Pshenichnyi, {\em Convex Analysis and Extremal Problems},  Nauka, Moscow, 1980. [In Russian]

\bibitem{Erm76}
Yu.M.~Ermoliev, {\em Methods of Stochastic Programming},  Nauka, Moscow, 1976. [In Russian]

\bibitem{Sho79}
N.Z.~Shor, {\em Minimization Methods for Non-Differentiable
Functions}, Naukova Dumka, Kiev, 1979
(Engl. transl. in Springer-Verlag, Berlin, 1985).

\bibitem{Pol83}
B.T.~Polyak, {\em Introduction to Optimization}, Nauka, Moscow, 1983
(Engl. transl. in Optimization Software, New York,  1987).

\bibitem{Pan85}
P.D.~Panagiotopoulos,
 {\em Inequality Problems in Mechanics and Their Applications},
 Birkhauser, Boston, 1985.

\bibitem{MN92}
M.M.~M\"{a}kela and P.~Neittaanm\"{a}ki,  {\em  Nonsmooth Optimization},
  World Scientific,  Singapore, 1992.

\bibitem{DV81}
V.F.~Dem'yanov and L.V.~Vasil'yev, {\em Nondifferentiable Optimization},
Nauka, Moscow, 1981 (Engl. transl. in Optimization Software,
New York, 1985).

\bibitem{ShS72}
N.Z.~Shor and L.P.~Shabashova, {\em Solution of minimax problems by the method
of generalized gradient descent with dilatation of the space}, Cybernetics, vol.8
(1972), pp.88--94.

\bibitem{NeS15}
Yu.~Nesterov and V.~Shikhman, {\em Quasi-monotone subgradient methods for nonsmooth
convex minimization},  J. Optim. Theory Appl., vol.165
(2015), pp. 917--940.

\end{thebibliography}
\end{document}